\documentclass[12pt,oneside,notitlepage]{article} 
\usepackage{amssymb,latexsym,amsmath,graphics,setspace}
\usepackage{epsfig,epsf,amsthm,amsfonts}

\newtheorem{corollary}{Corollary}
\newtheorem{definition}{Definition}
\newtheorem{example}{Example}
\newtheorem{lemma}{Lemma}

\newtheorem{theorem}{Theorem}
\begin{document}

\pagestyle{empty}
\begin{center}
\rule{165pt}{0pt} \\
\LARGE{A Short Study of Alexandroff Spaces}\\
\vspace{0.4cm}
\normalsize{by} \\
\large{Timothy Speer } \\
\vspace{5mm}
\normalsize{Department of Mathematics} \\
\normalsize{New York University} \\
\normalsize{tjs312@cims.nyu.edu} \\
\normalsize{July 2007} \\
\vspace{1.0cm}
\normalsize{\textbf{Abstract}}\\
\end{center}
In this paper, we discuss the basic properties of Alexandroff spaces. Several examples of Alexandroff spaces are given. We show how to construct new Alexandroff spaces from given ones. Finally, two invariants for compact Alexandroff spaces are defined and calculated for the given examples.  

\vspace{8.7cm}
\pagestyle{empty}             
\begin{center}
\copyright \quad Timothy Speer \\
All Rights Reserved, 2007
\end{center}
\pagebreak


\setcounter{page}{0} 
\tableofcontents

\pagebreak




\pagestyle{plain}
\section{Introduction}

The general definition of a topology is based off the properties of the standard Euclidean topology. The goal of this paper is to study spaces that have topologies which satisfy a stronger condition. Namely, arbitrary intersections of open sets are open. With this restriction, we lose important spaces such as Euclidean spaces, but the specialized spaces in turn display interesting properties that are not necessary for a standard topological space.

\begin{definition}
Let $X$ be a topological space, then $X$ is an Alexandroff space if arbitrary intersections of open sets are open. 
\end{definition}

Our first theorem shows that the condition of being an Alexandroff space severely limits the type of spaces we are considering. So we must look to more exotic examples of topologies to find useful spaces that are Alexandroff. A theorem proved later will allow us to construct Alexandroff spaces easily.   

\begin{lemma}
Any discrete topological space is an Alexandroff space. 
\end{lemma}
\begin{proof}
This is clear because in a discrete space any subset is open.  

\end{proof}

\begin{theorem}
Let $X$ be a metric space, then $X$ is an Alexandroff space iff $X$ has the discrete topology. 
\end{theorem}
\begin{proof}
$\Rightarrow$) Suppose $X$ is an Alexandroff space. Let $x$ be a point in $X$. Then the open balls $B(x,\frac{1}{n})$ with radius $\frac{1}{n}$ and center $x$, $n$ a natural number, are open in $X$. Since $X$ is an Alexandroff space, $\bigcap^{\infty}_{n=1}B(x,\frac{1}{n})$ is an open set. But by the properties of the metric, $\bigcap^{\infty}_{n=1}B(x,\frac{1}{n})=\{x\}$. So we have shown that singletons are open. Hence, $X$ has the discrete topology.

\vspace{5mm}
\noindent $\Leftarrow$) The reverse direction follows from Lemma 1.    

\end{proof}

The original definition given for an Alexandroff space is easy to state, however it is not too useful for proving theorems about Alexandroff spaces. To fix this we will use a different, yet equivalent definition. Namely, in an Alexandroff space, each point must have a minimal open set containing it. 

\begin{theorem}
$X$ is an Alexandroff space iff each point in $X$ has a minimal open neighborhood.
\end{theorem}
\begin{proof}
$\Rightarrow$) Suppose $X$ is an Alexandroff space with $x\in X$. Let
$O(x)=\{U\subset X: U$ is an open neighborhood of $x\}$. Take S(x)=$\bigcap U$ for
$U\in O(x)$, then $S(x)$ is an open neighborhood of $x$ because $X$ is Alexandroff.
And from the definition of $S(x)$, it is clear that $S(x)$ is a minimal open
neighborhood of $x$.

\vspace{5mm}
\noindent $\Leftarrow$) Suppose each $x\in X$ has a minimal open neighborhood S(x). Consider an
arbitrary intersection of open sets, $V=\bigcap_{\alpha\in A}U_\alpha$, where each
$U_\alpha$ is open in $X$. If $V=\emptyset$, then $V$ is open and we are done. If
$V\neq\emptyset$, then pick $x\in V$ and we have $x\in U_\alpha$ for all $\alpha\in
A$. Hence, $S(x)\subset U_\alpha$ for all $\alpha$ because $S(x)$ is the minimal
open neighborhood of $x$. Therefore, $S(x)\subset V$. Hence, $V$ is open because it
contains an open set around each of it's points.

\end{proof}

From here on, we will write $S(x)$ to denote the minimal open neighborhood of a
point $x$ in an Alexandroff space. 

\section{Minimal Open Neighborhoods}

In order to understand the properties of Alexandroff spaces, we will study the minimal open neighborhoods of these spaces. The proofs of the theorems show us that the minimal open neighborhoods of an Alexandroff space are the natural objects to study.

\begin{theorem} 
If $\beta$ is a collection of subsets of $X$ such that for each $x\in X$ there is a
minimal set $m(x)\in\beta$ containing $x$, then $\beta$ is a basis for a topology on
$X$ and $X$ is an Alexandroff space with this topology. In addition, $S(x)=m(x)$.
\end{theorem}
\begin{proof}
Clearly the sets in $\beta$ cover $X$. Suppose $U,V\in\beta$ and that $x\in U\cap
V$. Then $m(x)$ is a minimal set in $\beta$ containing $x$ so we must have
$m(x)\subset U$ and $m(x)\subset V$. Hence, $m(x)\subset U\cap V$. So $\beta$ is a
basis for a topology on $X$. To show that $X$ is an Alexandroff space with this
basis, take any $x\in X$ and let $U$ be any open set in $X$ containing $x$. Then
$U=\bigcup_{\alpha\in A}V_\alpha$ where $V_\alpha\in\beta$. But $x$ must be in
$V_\alpha$ for at least one $\alpha$, which means $m(x)\subset V_\alpha\subset U$.
Hence, m(x) is a minimal open set containing $x$. Therefore, $X$ is Alexandroff and
$S(x)=m(x)$.

\end{proof}

Theorem 3 is a powerful tool that will allow us to construct Alexandroff spaces by specifying a basis. 

\begin{theorem}If $X$ is an Alexandroff space with topology $T$, then $\beta=\{S(x):x\in X\}$ is a
basis for $T$.
\end{theorem}
\begin{proof}
By Theorem 3, $\beta$ is a basis for a topology $T'$ on $X$. We must show
$T'=T$. $T'\subset T$ because $\beta\subset T$. Suppose $U\in T$. Let
$U'=\bigcup_{x\in U}S(x)$, then $U'\in T'$ and we also have $U=U'$. Hence, $T\subset
T'$ so we have $T=T'$.

\end{proof}

\begin{corollary}
If $T$, $T'$ are two topologies on $X$ such that $X$ is an Alexandroff space and
$S_T(x)=S_{T'}(x)$ for all $x$ in $X$, then $T=T'$.
\end{corollary}
\begin{proof}
From Theorem 4, we have $\beta=\{S_T(x):x\in X\}$ and
$\beta'=\{S_{T'}(x):x\in X\}$ are bases for $T$ and $T'$ respectively. But
$\beta=\beta'$ so we have $T=T'$.     
 
\end{proof} 
 
We are now ready to construct some examples of Alexandroff spaces.

\begin{example} (An Alexandroff Topology on $\mathbb{R}^n$) 

\noindent Take $X$ to be $\mathbb{R}^n$ and let
$\beta=\{\overline{B(0,r)}:r\in\mathbb{R}_+\cup\{0\}\}$. Note that
$\overline{B(0,r)}$ is the closed ball with center $0$ and radius $r$ and that
$\overline{B(0,0)}=\{0\}$. If $x\in X$ then $\overline{B(0,|x|)}$ is a minimal set in $\beta$ containing $x$. $\beta$ satisfies the conditions of Theorem
3 so it is a basis for an Alexandroff topology on $X$.  Notice that this
space is not Hausdorff and not compact. 
\end{example}

\begin{example} (An Alexandroff Topology on $D^n$)

\noindent This topology is similar to the previous example on $\mathbb{R}^n$, except this time
we take $X=D^n$ and $\beta=\{\overline{B(0,r)}:r\in [0,1]\}$. With the topology
generated by $\beta$, $D^n$ is compact because any open cover of $D^n$ must contain
$D^n$. Hence, $D^n$ is a finite subcover. However, this space is also not Hausdorff.
\end{example}

\begin{example} (Disjoint Minimal Open Neighborhoods)

\noindent Take $X=\mathbb{R}\setminus\mathbb{Z}$ and $\beta=\{(n,n+1):n\in\mathbb{Z}\}$. Then
$X$ is an Alexandroff space with $S(x)=(n,n+1)$ where $n<x<n+1$. For any two minimal
open neighborhoods $S(x)\neq S(y)$ we have that $S(x)$ and $S(y)$ are disjoint.
\end{example}

\begin{example} (An Alexandroff Topology on $S^1$)

\noindent Take $X=S^1$ and define the sets $R_n=\{z\in S^1:z^n=1\}$ for
$n\in\mathbb{N}\cup\{0\}$. Let $\beta=\{R_n:n\in\mathbb{N}\cup\{0\}\}$. The point
$1\in S^1$ has the set $R_1=\{1\}$ as a minimal set in $\beta$ containing it.
Suppose $z\in S^1$ and $z\neq 1$. If $z^n=1$ holds only for $n=0$, then $R_0=S^1$ is
a minimal set in $\beta$ containing $z$. If $z^n=1$ for some $n\neq 0$ then let
$m=min\{n\in\mathbb{N}:z^n=1\}$. Suppose $z\in R_p$ for some $p>m$, then $p=qm+r$
where $0\leq r<m$. Then $1=z^p=(z^m)^qz^r=z^r$. But $r<m$ so it must be that $r=0$.
Hence, $p=qm$. If $\rho\in R_m$, then $\rho^p=(\rho^m)^q=1$ so $\rho\in R_p$. Hence,
$R_m\subset R_p$ and so $R_m$ is a minimal set in $\beta$ containing $z$. This space
is compact because any open cover contains $S^1$.      
\end{example}

We conclude this section with a simple statement about minimal open neighborhoods.

\begin{theorem}If $X$ is an Alexandroff space, then $S(x)$ is compact for all $x\in X$.
\end{theorem}
\begin{proof}
Let $\{V_\alpha\}_{\alpha\in A}$ be an open cover of $S(x)$. Then $x$ is in
$V_\alpha$ for some $\alpha\in A$. So we must have $S(x)\subset V_\alpha$. Hence,
$V_\alpha$ is a finite subcover of $\{V_\alpha\}_{\alpha\in A}$.

\end{proof}

\section{New Alexandroff Spaces From Old Spaces}

The tools in this section will allow us to construct new Alexandroff spaces from given ones. 

\begin{theorem}
If $X$ and $Y$ are Alexandroff spaces, then $X\times Y$ is also an Alexandroff space, with $S(x,y)=S(x)\times S(y)$. 
\end{theorem}
\begin{proof}
$X\times Y$ has as basis $\beta=\{U\times V: U$ is open in $X$ and $V$ is open in $Y\}$. Let $(x,y)\in X\times Y$, then $S(x)\times S(y)$ is in $\beta$ and the claim is that this is a minimal set in $\beta$ containing $(x,y)$. If $(x,y)\in U\times V\in\beta$, then $x\in U$ and $y\in V$ so $S(x)\subset U$ and $S(y)\subset V$. Therefore, $S(x)\times S(y)$ is contained in $U\times V$. So by applying Theorem 3, we know that $X\times Y$ is an Alexandroff space and $S(x,y)=S(x)\times S(y)$. 

\end{proof}

\begin{corollary}
If $X_1,...,X_n$ are Alexandroff spaces, then so is $X_1\times...\times X_n$. Furthermore, $S(x_1,...,x_n)=S(x_1)\times ...\times S(x_n)$.
\end{corollary}
\begin{proof}
Use induction and apply Theorem 6. 

\end{proof}

\begin{theorem}
If $A$ is a subspace of the Alexandroff space $X$, then $A$ is an Alexandroff space. In addition, $S_A(x)=A\cap S_X(x)$.  
\end{theorem}
\begin{proof}
Let $x\in A$ and suppose $U$ is an open neighborhood of $x$ in $A$. Then $U=A\cap V$ where $V$ is open in $X$. This means that $S_X(x)$ must be contained in $V$, so that $A\cap S_X(x)\subset A\cap V=U$. Hence, $A$ is an Alexandroff space with $S_A(x)=A\cap S_X(x)$.   

\end{proof}

\begin{theorem}
If $X/\sim$ is a quotient space of the Alexandroff space $X$, then $X/\sim$ is an Alexandroff space. 
\end{theorem}
\begin{proof}
Let $q:X\rightarrow X/\sim$ be the quotient map. Consider an arbitrary intersection, $\bigcap_{\alpha\in A}U_\alpha$, of open sets in $X/\sim$. We have $q^{-1}(\bigcap_{\alpha\in A}U_\alpha)=\bigcap_{\alpha\in A}q^{-1}(U_\alpha)$. Now, $q^{-1}(U_\alpha)$ is open in $X$ for each $\alpha\in A$ because q is the quotient map. Hence, $\bigcap_{\alpha\in A}q^{-1}(U_\alpha)$ is open in $X$ and therefore $\bigcap_{\alpha\in A}U_\alpha$ is open in $X/\sim$ by definition of the quotient topology. 
    
\end{proof}

\begin{example}
If $X$ is an Alexandroff space, then we can define an equivalence relation $\sim$ on $X$ by, $x\sim y$ iff $S(x)=S(y)$. We can then from the quotient space $X/\sim$, which is an Alexandroff space by Theorem 8. The following theorem will tell us exactly when $X/\sim$ is a discrete space, but we first need a definition. 
 
\end{example}

\begin{definition}
If $X$ is an Alexandroff space and $x\in X$, then $S(x)$ is called irreducible if $S(y)\subset S(x)$ implies $S(y)=S(x)$. 

\end{definition}

\begin{theorem}
$X/\sim$ is a discrete space iff $S(x)$ is irreducible for all $x$ in $X$. 
\end{theorem}
\begin{proof}
$\Rightarrow$) Suppose $X/\sim$ is discrete. If $q$ is the quotient map, then $q^{-1}([x])$ is open in $X$ for all $x$. Then $S(x)\subset q^{-1}([x])$ because $x$ is in $q^{-1}([x])$. If $y\in q^{-1}([x])$, then $x\sim y$ which means $S(x)=S(y)$. Therefore, $y$ is in $S(x)$ which gives us $q^{-1}([x])=S(x)$. Now suppose $S(z)\subset S(x)$, then $z$ is in $q^{-1}([x])$. So we must have $S(z)=S(x)$ because $z\sim x$. Hence, $S(x)$ is irreducible. 

\vspace{5mm}
\noindent$\Leftarrow$) Now suppose $S(x)$ is irreducible for all $x$ in $X$. Let $y\in q^{-1}([x])$, then $x\sim y$ so $S(x)=S(y)$. This gives us $q^{-1}([x])\subset S(x)$. Now if $y\in S(x)$, then $S(y)\subset S(x)$ and since $S(x)$ is irreducible, we must have $S(y)=S(x)$. Therefore, $x\sim y$ so that $y\in q^{-1}([x])$. So we have $q^{-1}([x])=S(x)$ is open in $X$. This means that $[x]$ is open in $X/\sim$ and so it is discrete.  

\end{proof}

\begin{theorem}
If $S(x)$ and $S(y)$ are distinct irreducible subsets of $X$, then $S(x)\cap S(y)=\emptyset$.
\end{theorem}
\begin{proof}
Suppose $S(x)$ and $S(y)$ are distinct irreducible subsets of $X$. Suppose that $z\in S(x)\cap S(y)$. Then $S(z)\subset S(x)\cap S(z)$ which implies that $S(z)\subset S(x)$ and $S(z)\subset S(y)$. Using the irreducibility of $S(x)$ and $S(y)$ gives us $S(x)=S(z)=S(y)$, which is a contradiction. Hence, $S(x)\cap S(y)=\emptyset$.

\end{proof}

\noindent\underline{Hausdorff Alexandroff Spaces}

\vspace{5mm}
For general topological spaces, the spaces that satisfy the Hausdorff property are the nicest to study. For Alexandroff spaces, the spaces that are Hausdorff are not very interesting to study as the following theorem shows.

\begin{theorem}
$X$ is a Hausdorff Alexandroff space iff for $x\neq y$ in $X$, we have $S(x)\cap S(y)=\emptyset$.
\end{theorem}
\begin{proof}
$\Rightarrow$) If $X$ is Hausdorff we can find disjoint open sets $U,V$ of $X$ such that $x\in U$ and $y\in V$. Then $S(x)\subset U$ and $S(y)\subset V$, so $S(x)$ and $S(y)$ must also be disjoint. 

\vspace{5mm}
\noindent $\Leftarrow$) This is trivial, just take $S(x)$ and $S(y)$ to be the disjoint open sets containing $x$ and $y$ respectively. 

\end{proof}

\begin{corollary}
$X$ is a Hausdorff Alexandroff space iff $X$ is discrete.
\end{corollary}
\begin{proof}
$\Rightarrow$) If $X$ is Hausdorff, then we claim that $S(x)=\{x\}$. To see this suppose $y\in S(x)$. Then $S(y)\subset S(x)$, which means $S(y)\cap S(x)=S(y)$. And since $S(y)\neq\emptyset$, we must have $y=x$ by Theorem 11. Hence, $\{x\}$ is open in $X$, so $X$ must be discrete. 

\vspace{5mm}
\noindent $\Leftarrow$). If $X$ is discrete, then it is clearly Hausdorff.  

\end{proof}

\section{Continuous Maps}

It seems reasonable to believe that being Alexandroff is a topological property of a space, and this is indeed the case. However, if $f:X\rightarrow Y$ is a continuous map and $X$ is an Alexandroff space, then $f(X)$ is not necessarily Alexandroff. For example, take $X=\mathbb{N}$ with the discrete topology and let $Y=\mathbb{Q}$ with the subspace topology from $\mathbb{R}$. Pick a bijection $f:\mathbb{N}\rightarrow\mathbb{Q}$, then $f$ must always be continuous but $f(\mathbb{N})=\mathbb{Q}$ is not an Alexandroff space. We need a stronger condition on $f$ to insure that $f(X)$ is an Alexandroff space.  

\begin{theorem}
Let $f:X\rightarrow Y$ be an open and continuous map. If $X$ is an Alexandroff space, then so is $f(X)$. In addition, if $y\in f(X)$, then $S(y)=f(S(x))$ where $f(x)=y$.  
\end{theorem}
\begin{proof}
Let $y\in f(X)$ and pick $x\in X$ such that $f(x)=y$. Then $f(S(x))$ is open in $f(X)$ because $f$ is an open map. Suppose that $y$ is contained in an open set $U$ in $f(X)$. This means $x\in f^{-1}(U)$ and since $f^{-1}(U)$ is open in $X$, we have $S(x)\subset f^{-1}(U)$. Therefore, $f(S(x))\subset U$, which shows $f(X)$ is Alexandroff with $S(y)=f(S(x))$. 

\end{proof} 

\begin{corollary}
If $X$ is homeomorphic to $Y$ and $X$ is an Alexandroff space, then so is $Y$. 
\end{corollary}
\begin{proof}
If $f$ is a homeomorphism between $X$ and $Y$, then $f$ is open and continuous with $f(X)=Y$. By Theorem 12, $Y$ is an Alexandroff space.   

\end{proof}

\begin{theorem}
If $X$ and $Y$ are Alexandroff spaces that satisfy the following:

\vspace{3mm}
i.) There is a bijection $b:\{S(x):x\in X\}\rightarrow\{S(y):y\in Y\}$

ii.) There are homeomorphisms $f_x:S(x)\rightarrow b(S(x))$ for each $x$ in $X$.

iii.) If $S(x_1)\cap S(x_2)\neq\emptyset$ then $f_{x_1}(S(x_1)\cap S(x_2))=f_{x_2}(S(x_1)\cap S(x_2))$.

\vspace{3mm}
\noindent then $X$ and $Y$ are homeomorphic. 
 
\end{theorem}
\begin{proof}
Define $h:X\rightarrow Y$ by $h(x)=f_z(x)$ where $x\in S(z)$. We must show that $h$ is well-defined. If $x\in S(z_1)$ and $x\in S(z_2)$, then $f_{z_1}(x)=f_{z_2}(x)$ by $iii.)$. So the choice for $h(x)$ does not depend on which $z$ we choose. Now suppose $y\in Y$, then there exists $x\in X$ such that $b(S(x))=S(y)$. Then $f_x:S(x)\rightarrow S(y)$ is a homeomorphism so there exists $z\in S(x)$ such that $f_x(z)=y$. So $h(z)=f_x(z)=y$, which means $h$ is onto. Suppose $h(x_1)=h(x_2)=y$, then $y\in b(S(x_1))$ so we must have $S(y)\subset b(S(x_1))$. But since $f^{-1}_{x_1}(S(y))$ is open in $S(x_1)$ and $x_1\in f^{-1}_{x_1}(S(y))$, we must have $S(x_1)=f^{-1}_{x_1}(S(y))$. Therefore, $b(S(x_1))=S(y)$. Similarly, we get that $b(S(x_2))=S(y)$. And $b$ is injective so $S(x_1)=S(x_2)$. This tells us that $f_{x_1}=f_{x_2}$, so that $f_{x_1}(x_1)=f_{x_1}(x_2)$. But this map is injective so we must have $x_1=x_2$. So $h$ is also injective. To show continuity, we only need to use the basis of $Y$ of minimal open neighborhoods. $h^{-1}(S(y))=S(x)$, where $h(x)=y$ so this is an open set. Hence, $h$ is continuous. Similarly, $h^{-1}$ is continuous.          

\end{proof}

\section{Compact Alexandroff Spaces} 

We will begin the study of compact Alexandroff spaces. If $X$ is a compact Alexandroff space, then it is covered by the set $\{S(x):x\in X\}$. Hence, it must be covered by a finite number of minimal open neighborhoods. This key property will allow us to define several invariants of compact Alexandroff spaces.

\begin{definition}
If $X$ is a compact Alexandroff space, then we define $min(X)=min\{|V|:V$ is a finite cover of $X$ by minimal open neighborhoods$\}$. 
\end{definition}  

\begin{example}
The topologies given in Example 2 and Example 4 on $D^n$ and $S^1$ are both compact. We have that $min(D^n)=1$ and $min(S^1)=1$.  
\end{example}

The following theorem shows us that the $min(X)$ of an Alexandroff space, is an invariant of that space. However, it is not enough to distinguish between the two spaces given in examples 2 and 4.    

\begin{theorem}
If $X$ and $Y$ are compact Alexandroff spaces that are homeomorphic, then $min(X)=min(Y)$. 
\end{theorem}
\begin{proof}
Let $h:X\rightarrow Y$ be a homeomorphism. Let $\{S(y_1),...,S(y_{min(Y)})\}$ be an open cover of $Y$ by minimal open neighborhoods. Then $\{S(x_1)=h^{-1}(S(y_1)),...,S(x_{min(Y)})=h^{-1}(S(y_{min(Y)})\}$, where $h(x_i)=y_i$, is an open cover of $X$ by minimal open neighborhoods. Therefore, $min(X)\leq min(Y)$. By a similar argument in the other direction, we get $min(Y)\leq min(X)$. Hence, $min(X)=min(Y)$.    

\end{proof}

As seen in Theorem 13, we can define a homeomorphism between two Alexandroff spaces by specifying homeomorphisms between the minimal open neighborhoods of the spaces. This shows that minimal open neighborhoods are the natural objects of study in Alexandroff spaces. For compact Alexandroff spaces, there is an even more basic object of study. 

\begin{definition}
If $X$ is an Alexandroff space, then $S(x)$ is called basic if for any $S(y)$; if $S(x)\subset S(y)$ and $S(z)\subset S(y)$ then $S(x)\subset S(z)$ and if $S(x)$ is not contained in $S(y)$, then $S(x)$ and $S(y)$ are disjoint. 
\end{definition} 

\begin{example}
In Example 1 and Example 2, the only basic set is the set $\{0\}$. In Example 4, the only basic set is $\{1\}$. Example 3 has an infinite number of basic sets because each minimal open neighborhood is also a basic set.  
\end{example}

\begin{theorem}
If $S(x)$ is a basic subset of $X$, then $S(x)$ is irreducible. 
\end{theorem}
\begin{proof}
Suppose that $S(y)\subset S(x)$. Since $S(x)$ is basic, if $S(x)$ is not contained in $S(y)$ we must have $S(x)\cap S(y)=\emptyset$. This can't be because $y$ is in $S(y)$ and $S(x)$. Hence, $S(x)\subset S(y)$, which implies $S(y)=S(x)$.  

\end{proof}

\begin{corollary}
If $S(x)$ and $S(y)$ are basic subsets of $X$, then they are disjoint.
\end{corollary}
\begin{proof}
First apply Theorem 15 and then Theorem 10. 

\end{proof}

\begin{theorem}
If $X$ is an Alexandroff space, then $S(x)$ contains at most one basic set for each $x\in X$. 
\end{theorem}
\begin{proof}
Suppose $S(x)$ contains two basic sets $S(y)$ and $S(z)$. Then by definition of basic, $S(y)\subset S(z)$ and $S(z)\subset S(y)$. Hence, $S(y)=S(z)$. 

\end{proof}

The number of basic sets in a compact Alexandroff space will give us another invariant of the space. We must first show that there is always a finite number of basic sets in a compact Alexandroff space. 

\begin{theorem}
If $X$ is a compact Alexandroff space, then the number of basic sets is less than or equal to $min(X)$. 
\end{theorem}
\begin{proof}
Cover $X$ by $\{S(x_1),...,S(x_{min(X)})\}$. Suppose the number of basic subsets of $X$ is greater than $min(X)$. We claim that each basic set, $S(y)$, is contained in $S(x_i)$ for some $i$. Suppose not, then $S(y)\cap S(x_i)=\emptyset$ for all $i$. This is a contradiction because the $S(x_i)$ cover $X$. Hence, $S(y)\subset S(x_i)$ for some $i$. If the number of basic subsets of $X$ is greater than $min(X)$, then some $S(x_i)$ contains more than one basic set which contradicts Theorem 16. So we must have the number of basic subsets of $X$ less than or equal to $min(X)$.  

\end{proof}

\begin{definition}
If $X$ is a compact Alexandroff space, then define $index(X)$ to be the number of basic subsets of $X$. 
\end{definition}

\begin{example}
Then index of the spaces in examples 2 and 4 is 1. 
\end{example}

\begin{theorem}
If $X$ and $Y$ are compact Alexandroff spaces that are homeomorphic, then $index(X)=index(Y)$. 
\end{theorem}
\begin{proof}
Let $h:X\rightarrow Y$ be a homeomorphism. Suppose $S(x)$ is a basic subset of $X$. Then we claim that $h(S(x))=S(h(x))$ is a basic subset of $Y$. Suppose $S(h(x))\subset S(y)$ and $S(z)\subset S(y)$, then $S(h^{-1}(y))=h^{-1}(S(y))\supset h^{-1}(S(z))=S(h^{-1}(z))$ and $S(x)\subset S(h^{-1}(y))$. So we must have $S(x)\subset S(h^{-1}(z))$ which means $S(h(x))\subset S(z)$. If $S(y)$ does not contain $S(h(x))$, then $S(h^{-1}(y))$ does not contain $S(x)$, which implies $S(h^{-1}(y))\cap S(x)=\emptyset$. Therefore, $S(y)\cap S(h(x))=\emptyset$. So $S(h(x))$ is a basic subset of $Y$. This gives us $index(X)\leq index(Y)$. By a similar argument, we obtain $index(Y)\leq index(X)$. So $index(X)=index(Y)$.  

\end{proof}

\include{references} 
\end{document}